\input amstex.tex
\documentstyle{amsppt}
\magnification1200
\hsize=12.5cm
\vsize=18cm
\hoffset=1cm
\voffset=2cm

\def\DJ{\leavevmode\setbox0=\hbox{D}\kern0pt\rlap
{\kern.04em\raise.188\ht0\hbox{-}}D}
\def\dj{\leavevmode
 \setbox0=\hbox{d}\kern0pt\rlap{\kern.215em\raise.46\ht0\hbox{-}}d}

\def\txt#1{{\textstyle{#1}}}
\baselineskip=13pt
\def\hf{{\textstyle{1\over2}}}

\def\d{{\,\roman d}}
\def\e{\varepsilon}
\def\D{\Delta}

\def\f{\varphi}
\def\G{\Gamma}
\def\k{\kappa}
\def\s{\sigma}

\def\={\;=\;}
\def\zx{\zeta(\hf+ix)}
\def\zt{\zeta(\hf+it)}

\def\D{\Delta}

\def\R{\Re{\roman e}\,} \def\I{\Im{\roman m}\,}
\def\z{\zeta}

\def\hf{{\textstyle{1\over2}}}
\def\txt#1{{\textstyle{#1}}}
\def\f{\varphi}
\def\Z{{\Cal Z}}
\def\M{{\Cal M}}
\def\le{\leqslant}
\def\ge{\geqslant}
\font\tenmsb=msbm10
\font\sevenmsb=msbm7
\font\fivemsb=msbm5
\newfam\msbfam
\textfont\msbfam=\tenmsb
\scriptfont\msbfam=\sevenmsb
\scriptscriptfont\msbfam=\fivemsb
\def\Bbb#1{{\fam\msbfam #1}}

\def \NN {\Bbb N}
\def \CC {\Bbb C}
\def \RR {\Bbb R}
\def \ZZ {\Bbb Z}

\font\ff=cmr8
\def\txt#1{{\textstyle{#1}}}
\baselineskip=13pt

\font\teneufm=eufm10
\font\seveneufm=eufm7
\font\fiveeufm=eufm5
\newfam\eufmfam
\textfont\eufmfam=\teneufm
\scriptfont\eufmfam=\seveneufm
\scriptscriptfont\eufmfam=\fiveeufm
\def\mathfrak#1{{\fam\eufmfam\relax#1}}

\font\tenmsb=msbm10
\font\sevenmsb=msbm7
\font\fivemsb=msbm5
\newfam\msbfam
     \textfont\msbfam=\tenmsb
      \scriptfont\msbfam=\sevenmsb
      \scriptscriptfont\msbfam=\fivemsb
\def\Bbb#1{{\fam\msbfam #1}}

\def \NN {\Bbb N}
\def \CC {\Bbb C}

\def \RR {\Bbb R}
\def \ZZ {\Bbb Z}

  \def\rightheadline{{\hfil{\ff
Mellin transforms of powers of Hardy's function}\hfil\tenrm\folio}}

  \def\leftheadline{{\tenrm\folio\hfil{\ff
 Aleksandar Ivi\'c }\hfil}}
  \def\emptyheadline{\hfil}
  \headline{\ifnum\pageno=1 \emptyheadline\else
  \ifodd\pageno \rightheadline \else \leftheadline\fi\fi}

\topmatter
\title
On the Mellin transforms of powers of Hardy's function
\endtitle
\author
 Aleksandar Ivi\'c
\endauthor
\address
Katedra Matematike RGF-a,
Universitet u Beogradu,  \DJ u\v sina 7,
11000 Beograd, Serbia.\bigskip
\endaddress
\subjclass 11 M 06
\endsubjclass
\email
{\tt
ivic\@rgf.bg.ac.rs,\enskip aivic\@matf.bg.ac.rs}
\endemail
\dedicatory
\enddedicatory
\abstract
Various properties of the  Mellin transform function
$$
{\Cal M}_k(s) :=
\int_1^\infty Z^k(x)x^{-s}\d x
$$
are investigated, where
$$
Z(t) := \zt{\bigl(\chi(\hf+it)\bigr)}^{-1/2}, \quad \z(s) = \chi(s)\z(1-s)
$$
is Hardy's function. Connections with power moments of $|\zt|$
are established, and natural boundaries of ${\Cal M}_k(s)$ are discussed.
\endabstract
\endtopmatter
\document
\head
1. Introduction
\endhead
Power moments of $|\zt|$ are a central problem in the theory of the
Riemann zeta-function $$
\z(s) = \sum\limits_{n=1}^\infty n^{-s}\qquad(\s = \R s >1),
$$
which admits analytic continuation to $\CC$, having only a simple pole at $s=1$.
A vast literature
exists on this subject (see e.g., the monographs [10], [11], [26] and [28]).
One way to tackle them is to deal with the (modified) Mellin transform function
$$
\Z_k(s) := \int_1^\infty |\zx|^{2k}x^{-s}\d x\qquad(k\in\NN),\leqno(1.1)
$$
where $\s = \R s$ is so large that the integral in (1.1) converges
absolutely. These functions in the cases when $k=1$ or $k=2$ have been
intensively investigated (e.g., see the works [14], [18], [23] and [24]).
It is known that $\Z_1(s)$ has meromorphic continuation to $\CC$. It has a
pole of order two at $s=1$ and the principal part of the Laurent expansion
at $s=1$ of $\Z_1(s)$ is
$$
{1\over(s-1)^2} + {2\gamma-\log(2\pi)\over s-1},
$$
where $\gamma = -\G'(1) = 0.577215\ldots\,$ is Euler's constant. It also has
simple poles at $s =-1,-3,\ldots\,$, whose residues can be expressed
explicitly in term of Bernoulli numbers (see M. Lukkarinen [23]).

\medskip
The analytic continuation of $\Z_2(s)$ has also
(see e.g., Y. Motohashi [24]) infinitely many poles. Namely
in the half-plane $\R s > 0$ it has
the following singularities: the pole $s = 1$ of order five,  simple
poles at $s = {1\over2} \pm i\k_j\,\left(\k_j =\sqrt{\lambda_j -
{1\over4}}\,\right)$ and poles at $s =
\rho/2$.  Here $\rho$ denotes complex zeros of
$\zeta(s)$, and $\,\{\lambda_j = \k_j^2 +
{1\over4}\} \,\cup\, \{0\}\,$ is the discrete spectrum of the
non-Euclidean Laplacian acting on $SL(2,\ZZ)$-automorphic forms.
This shows that $\Z_2(s)$ has a different and more complicated
structure than $\Z_1(s)$.

\medskip
Instead of $\Z_k(s)$  one can consider the more general
Mellin transform function
$$
{\Cal M}_k(s) :=
\int_1^\infty Z^k(x)x^{-s}\d x\qquad(k\in\NN),
\leqno(1.2)
$$
where again $\s = \R s$ is so large that the integral in (1.2) converges
absolutely. Here $Z(x)$ is the classical Hardy function, defined as
$$
Z(t) := \zt\bigl(\chi(\hf+it)\bigr)^{-1/2}, \quad \z(s) = \chi(s)\z(1-s),\leqno(1.3)
$$
with
$$\chi(s) = 2^s\pi^{s-1}\sin(\hf \pi s)\G(1-s),\quad
\chi(s)\chi(1-s)=1.
$$
It follows  that $\overline{\chi(\hf + it)} = \chi(\hf-it)$, so that
$Z(t)\in\RR$ when $t\in\RR$
and $|Z(t)| =|\zt|$. Thus the zeros of $\z(s)$ on the ``critical line'' $\R s =1/2$
correspond to the real zeros of $Z(t)$, which makes $Z(t)$ an invaluable tool
in the study of the zeros of the zeta-function on the critical line.
Note that when $k = 2\ell$ is even, then
$$
{\Cal M}_{2\ell}(s) = \int_1^\infty Z^{2\ell}(x)x^{-s}\d x
= \int_1^\infty |\zx|^{2\ell}x^{-s}\d x = {\Cal Z}_\ell(s)
$$
in former notation. Hence ${\Cal M}_k(s)$ is also closely connected to the moments
of $\zx$. If we define, for fixed $k\in\NN$, the $k$--th moment
of Hardy's function as
$$
{\Cal I}_k(x) := \int_1^xZ^k(y)\d y,\leqno(1.4)
$$
then on integrating by parts we find that
$$
{\Cal M}_k(s) = s\int_1^\infty {\Cal I}_k(x)x^{-s-1}\d x,
\leqno(1.5)
$$
so that the properties of ${\Cal I}_k(x)$ are reflected on ${\Cal M}_k(s)$.
Conversely, the Mellin inversion formula  gives
$$
Z^k(x) = {1\over2\pi i}\int_{(c)}{\Cal M}_k(s)x^{s-1}\d s
\leqno(1.6)
$$
for suitable $c\;(>0)$. From (1.6) we obtain by integration
$$
{\Cal I}_k(x) = {1\over2\pi i}\int_{(c)}{\Cal M}_k(s){x^s\over s}\d s
+ O(1).
$$

\smallskip
The plan of the paper is as follows. In Section 2 we consider $\M_k(s)$
and obtain some general results for this modified Mellin transform. Section
3 is devoted to $\M_k(s)$ in the special cases when $k=1$ and $k=3$.
The discussion related to the analytic continuation of $\M_3(s)$ is made
by the use of the cubic moment of $Z(x)$, which is dealt with in Section 4.
Finally the natural boundaries of $\M_k(s)$ and related problems
are treated in Section 5.
\head
2. Properties of $\M_k(s)$
\endhead

First we recall that the Mellin transform of $f(x)$ is commonly defined as
$$
\M[f(x)] =  F(s) := \int_0^\infty f(x)x^{s-1}\d x \qquad(s = \s + it).\leqno(2.1)
$$
Mellin and Laplace transforms play an important r\^ole in Analytic Number Theory.
They can be viewed, by  a change of variable, as special cases of
Fourier transforms, and their properties can be deduced from the
general theory of Fourier transforms (see e.g., E.C. Titchmarsh [27]).
For example, by the change of variable $x = {\roman e}^t,\,z = s-1$, (1.2) becomes
$$
\int_0^\infty Z^k({\roman e}^t){\roman e}^{-zt}\d t\qquad
(\R z > 0),
$$
which is the Laplace transform of $Z^k({\roman e}^t)$.
The reason that we have defined in (1.1) and (1.2) somewhat differently
the Mellin transforms $\Z_k(s), \M_k(s)$ is practical: the lower limit of integration
$x=1$ dispenses with potential convergence problems at $x=0$,
while the appearance
of $x^{-s}$ instead of the familiar $x^{s-1}$ stresses the analogy with Dirichlet
series where one has a sum of $f(n)n^{-s}$ and not $f(n)n^{s-1}$.

\smallskip
One of the basic properties of Mellin transforms is the  inversion formula
$$
\hf\{f(x+0) + f(x-0)\} = {1\over2\pi i}\int\limits_{(\s)}F(s)x^{-s}\d s =
{1\over2\pi i}\lim_{T\to\infty}\int\limits_{\s-iT}^{\s+iT}F(s)x^{-s}\d s.\leqno(2.2)
$$
Formula (2.2) certainly holds  if $f(x)x^{\s-1} \in L(0,\infty)$, and $f(x)$
is of bounded variation on every finite $x$--interval. Therefore
the inversion formula (1.6) follows from (2.2) by a change of variable.
Note that if $G(s)$ denotes the Mellin transform of $g(x)$  then,
assuming $f(x)$ and $g(x)$ to be real-valued, we formally have
$$
\eqalign{& {1\over2 \pi i}\int_{(\s)}F(s)\overline{G(s)}\d s
=\int_0^\infty g(x)\left({1\over2 \pi
i}\int_{(\s)}F(s)x^{\s-it-1}\d s\right)\d x\cr  &= \int_0^\infty
g(x)x^{2\s-1}\left({1\over2 \pi i}\int_{(\s)} F(s)x^{-s}\d s\right)\d x
=\int_0^\infty f(x)g(x)x^{2\s-1}\d x.\cr}
\leqno(2.3)
$$
The relation (2.3) is a form of Parseval's formula for Mellin
transforms, and it offers various possibilities  for mean square
bounds. A condition under which (2.3) holds is that $x^\s f(x)$ and $x^\s g(x)$
belong to $L^2((0,\infty),\, \d x/x)$. A variant of (2.3) is (see [27, Theorem 73])
$$
{1\over2\pi i}\int_{(c)}F(w)G(s - w)\d w    \=
\int_0^\infty f(x)g(x)x^{s-1}\d x,\leqno(2.4)
$$
which holds if $x^cf(x)$ and $x^{\s-c}g(x)$ belong to $L^2((0,\infty),\,
\d x/x)$, where as usual
$$
L^p(a,b) := \left\{\,f(x) \;\Biggl|\; \int_a^b|f(x)|^p\d x < \infty\;\right\}.
$$
Our first result is

\smallskip
THEOREM 1. {\it For $c \ge c_k >0, k\ge 2$ and $\s = \R s \ge \s_1(k)\;(>1)$ we have
$$
\M_k(s) = {1\over2\pi i}\int_{(c)}\M_{k-r}(w)\M_r(1-w+s)\d w \quad(r = 1,\ldots,k-1).
\leqno(2.5)
$$
In particular, for $\s > c > 1$,}
$$
\M_3(s)= {1\over2\pi i}\int_{(c)}\M_1(w)\M_2(1-w+s)\d w.\leqno(2.6)
$$

\medskip
{\bf Proof}. Consider
$$
f(x) = Z^{k-r}\Bigl({1\over x}\Bigr){1\over x},
\quad g(x) = Z^{r}\Bigl({1\over x}\Bigr){1\over x}\quad\quad(0<x\le 1),
$$
and $f(x) = g(x) =0$ if $x>1$. With the change of variable $y = 1/x$ we have
$$
F(s)= \int\limits_0^\infty f(x)x^{s-1}\d x=
\int\limits_0^1Z^{k-r}\Bigl({1\over x}\Bigr)x^{s-2}\d x = \int\limits_1^\infty
Z^{k-r}(y)y^{-s}\d y = \M_{k-r}(s),
$$
and likewise $G(s) = \M_r(s)$. Hence (2.4) yields
$$
{1\over2\pi i}\int_{(c)}\M_{k-r}(w)\M_r(s-w)\d w =
\int_0^1Z^{k-r}\Bigl({1\over x}\Bigr)Z^{r}\Bigl({1\over x}\Bigr)x^{s-3}\d x
= \M_k(s-1),
$$
again with the change of variable $y = 1/x$. Finally changing $s-1$ to $s$ we
obtain (2.5).

\smallskip
To establish (2.6) let $I$ denote the integral on the right-hand side.
We shall use the following elementary (see [12, Lemma 4])

\smallskip
LEMMA 1. {\it Suppose that $g(x)$ is a real-valued,
integrable function on $[a,b]$, a subinterval
of $[1,\,\infty)$, which is not necessarily finite. Then}
$$
\int\limits_0^{T}\left|\int\limits_a^b g(x)x^{-s}\d x\right|^2\d t
\le 2\pi\int\limits_a^b g^2(x)x^{1-2\s}\d x \quad(s = \s+it\,,T > 0,\,a<b).
\leqno(2.7)
$$

\smallskip
Then by using the Cauchy-Schwarz inequality, the well-known bounds (see [10])
$$
\int_0^T |\zt|^{2k}\d t \;\ll\; T(\log T)^{k^2}\qquad(k = 1,2)\leqno(2.8)
$$
and (2.7) (considering $t = \I s$ fixed and letting $T\to\infty$)  we obtain
$$
\eqalign{
I^2 & \leq \int_{-\infty}^\infty |\M_1(c+iv)|^2\d v
\int_{-\infty}^\infty |\M_1(1-c+\s+i(v+t))|^2\d v\cr&
\ll \int_1^\infty |\zx|^2x^{1-2c}\d x
\int_1^\infty |\zx|^4x^{2c-2\s-1}\d x\cr&
\ll 1,\cr}
$$
since $1-2c< -1, 2c -2\s - 1 < -1$. Therefore $I$ converges absolutely
and (2.6) holds, providing incidentally the analytic continuation of
$\M_3(s)$ to $\s>1$ (this also follows directly from the defining relation (1.2)).

\medskip
 THEOREM 2. {\it If $k = 1,2,3,4$ and $c>1$ is fixed, then for $U\gg x$
 and $\e>0$ sufficiently small we have}
 $$
 Z^k(x) = {1\over2\pi i}\int_{c-iU}^{c+iU}x^{s-1}\M_k(s)\d s + O_{\e,k}(x^{c-1}U^{-\e/2}).
 \leqno(2.9)
 $$

 \medskip
 {\bf Proof}. In view of (2.8) $\M_k(s)\;(k\le4)$ converges absolutely for $\s>1$.
 Hence the inversion formula (1.6) and the residue theorem yield
$$
\eqalign{
Z^k(x) &= {1\over2\pi }\int_{(c)}x^{s-1}\M_k(s)\d s\cr&
=  {1\over2\pi }\left(\int_{c-iU}^{c+iU} + \int_{c-i\infty}^{c-iU}
+ \int_{c+iU}^{c+i\infty}\right) + O_{\e,k}(x^{\e})\cr&
= {1\over2\pi }(I_1 + I_2+I_3) + O_{\e,k}(x^\e),\cr}
$$
say. Here and later $\e \,(>0)$ denotes constants which may be arbitrarily small,
but are not necessarily the same ones at each occurrence.
The $O$-term comes from the residue at $s=1$ (for $k=1$ the function
$\M_1(s)$ is regular for $s=1$, while for $k=3$ very
likely $\M_3(s)$ is also regular at $s=1$, but this has not been proved yet).
Therefore to prove (2.9) it suffices to show that
$$
I_3 \;\ll_{\e,k}\; x^{c-1}U^{-\e/2},\leqno(2.10)
$$
since the estimation of $I_2$ is analogous to the estimation of $I_3$.
For $\s>1, T_1\le t\le 2T_1 $ (with the aim of taking later $T_1 = U, T_1 = 2U$ etc.)
we have
$$
\M_k(s) = \int_1^{T_1^{1-\e}}Z^k(u)\f(u)u^{-s}\d u +
\int_{{1\over2}T_1^{1-\e}}^\infty
Z^k(u)(1-\f(u))u^{-s}\d u = I_4+I_5,
$$
say. Here $\f(u)\;(\ge0)$ is a smooth function supported in $[1, T_1^{1-\e}]$ such
that $\f(u) = 1$ for $1\le u \le \hf T_1^{1-\e}$ and
$$
\f^{(r)}(u) \;\ll_r\; T_1^{r(1-\e)} \qquad(r=0,1,2,\ldots\,).\leqno(2.11)
$$

Repeated integration by parts shows that, for $N\ge N_0(\e,k)$,
$$
\eqalign{
I_4 &= {c_{1,k}\over s-1} + {1\over s-1}\int_1^{T_1^{1-\e}}u^{1-s}
\left(\f(u)Z^k(u)\right)^{'}\d u = \ldots \cr&
= {c_{1,k}\over s-1}+ \ldots + {c_{N,k}\over (s-1)^N} + O_{N,k}(T_1^{-{1\over2}\e N})
\cr}\leqno(2.12)
$$
since, for $\ell_j, m_j\ge0$ and $\ell_1 + \ldots + \ell_N = k$ (for a formula for
$Z^{(m)}(x)$ see [21, p. 87]; see also [18, p. 313]),
$$
\int_1^X \Bigl(Z^{\ell_1}(x)\Bigr)^{(m_1)}
\ldots \Bigl(Z^{\ell_N}(x)\Bigr)^{(m_N)}\d x \ll_{\e,k,N} X^{1+\e}.\leqno(2.13)
$$
One obtains (2.13) similarly as (2.8), using H\"older's inequality,
the defining relation (1.3)
and the asymptotics of the $\chi$--function. The reason that we do not have
(yet) Theorem 2 for $k>4$ is essentially the fact that we do not have yet the bound
$$
\int_0^T|\zt|^m\d t \;\ll_\e\; T^{1+\e}
$$
for any fixed $m>4$.

Hence by the first derivative test
$$
\int_{c+iT_1}^{c+i2T_1}x^{s-1}I_4\d s = {c_{1,k}\over2\pi i}\int_{c+iT_1}^{c+i2T_1}
{x^{s-1}\over s-1}\d s + O\left({x^{c-1}\over T_1}\right) =
O\left({x^{c-1}\over T_1}\right).
$$
On the other hand
$$
\eqalign{
\int_{c+iT_1}^{c+i2T_1}x^{s-1}I_5\d s &
= i\int_{T_1}^{2T_1}x^{c+it-1}\left(\int_{{1\over2}T_1^{1-\e}}^\infty
Z^k(u)(1-\f(u))u^{-c-it}\d u\right)\d t\cr &
= ix^{c-1}\int_{{1\over2}T^{1-\e}}^\infty Z^k(u)(1-\f(u))u^{-c}
\left(\int_{T_1}^{2T_1}{\roman e}^{it\log(x/u)}\d t\right)\d u.
\cr}
$$
For $T_1\gg x$ it follows,  by direct integration, that the last integral
over $t$ is bounded. Thus the last expression, for some constant $c_k\ge0$,
is
$$
\ll \; x^{c-1}(\log T_1)^{c_k}T_1^{(1-\e)(1-c)}.
$$
Therefore we have
$$
I_3 \ll {x^{c-1}\over U} + x^{c-1}(\log U)^{c_k}U^{(1-\e)(1-c)} \ll x^{c-1}U^{-\e/2}
$$
if $\e>0$ is sufficiently small, and (2.9) follows. Theorem 2 is proved.

\medskip
{\bf Remark 1}. We can get (at least in principle) the information about
the sixth moment of $\zt$ from ${\Cal M}_3(s)$.
Namely from (2.9) with $k=3$ or from the method of proof of Lemma 4 of [10]
 we get that
$$
\int_T^{2T}|\zt|^6\d t \ll_\e T^{2\s-1}
\int_1^{T^{1+\e}}|{\Cal M}_3(\s+it)|^2\d t
+ T^{1+\e}\quad(\hf < \s \le 1), \leqno(2.14)
$$
provided that  ${\Cal M}_3(s)$ can be continued to $\R s \ge\s$
(and that is the catch!).
Heuristically, we should be able to have $\s = 3/4+\e$, and then the integral on the
right-hand side of (2.14) should be $\ll_\e T^{1/2+\e}$,
giving a weak form of the sixth moment.
Note that  (see [12, eq. (4.7)]) for the eighth moment we have
$$
\int_T^{2T}|\zt|^8\d t \ll_\e
T^{2\s-1}\int_1^{T^{1+\e}}|\Z_2(\s+it)|^2\d t + T^{1+\e}\quad(\hf < \s \le1),
\leqno(2.15)
$$
and an analogue of (2.14) and (2.15) holds also for the mean square and
fourth power of $|\zt|$. In these cases, however, the results are not
of particular interest, since we have precise information which has been
obtained by other methods. The bounds for the sixth moment of $|\zt|$ are
intricately connected with the problem of the analytic continuation of
$\M_3(s)$ to the region $\s \le 1$. It should be noted that the bounds
$$
\int_0^T|\zt|^8\d t \;\ll_{\e}\; T^{1+\e}
$$
and
$$
\eqalign{
\int_T^{2T}|\Z_2(\s + it)|^2\d t &\;\ll_\e\;
T^{4-4\s+\e}\qquad(\hf < \s \le 1),\cr
\int_T^{2T}|\Z_2(\s + it)|^2\d t &\;\ll_\e\;
T^{2-2\s+\e} + T^{-1}\quad(\s \ge 1).
\cr}
$$
are equivalent (see [12, eqs. (4.3) and (4.8)]).

\medskip
The next result is a generalization of Theorem 4 of [15]. This is

\medskip
THEOREM 3. {\it In the region of absolute convergence we have}
$$
\M_k^2(s) = 2\int_1^\infty x^{-s}\left(\int_{\sqrt{x}}^x Z^k(u)
Z^k\Bigl({x\over u}\Bigr){\d u\over u}\right)\d x.
\leqno(2.16)
$$

\bigskip
{\bf Proof of Theorem 3}. Set $f(x) = Z^k(x)$ and make the change
of variables $xy = X,\,x/y = Y$, so that the absolute value of the
Jacobian of the transformation is equal to $1/(2Y)$. Therefore
$$\eqalign{
\M_k^2(s) &= \int_1^\infty\int_1^\infty (xy)^{-s}f(x)f(y)\d x\d y\cr&
= {1\over2}\int_1^\infty X^{-s}\int_{1/X}^X{1\over Y}
f(\sqrt{XY}\,)f(\sqrt{X/Y}\,)\d Y\d X.\cr}
$$
But as we have ($y = 1/u$)
$$
\int_{1/x}^x
f(\sqrt{xy}\,)f(\sqrt{x/y}\,){\d y\over y} = \int_{1/x}^1 + \int_1^x =
2\int_1^xf(\sqrt{x/u})f(\sqrt{xu}\,){\d u\over u},
$$
we obtain that, in the region of absolute convergence, the identity
$$
\M_k^2(s) = \int_1^\infty x^{-s}\left(\int_1^x f(\sqrt{xy}\,)
f(\sqrt{x/y}\,){\d y\over y}\right)\d x
$$
is valid. The inner integral here becomes, after the change of
variable $\sqrt{xy} = u$,
$$
2\int_{\sqrt{x}}^x f(u)f\bigl({x\over u}\bigr){\d u\over u},
$$
and (2.16) follows. The argument also shows that, for  $0 < a < b$
and any integrable function $f$ on $[a,\,b]$,
$$
\left(\int_a^b f(x)x^{-s}\d x\right)^2
= 2\int_{a^2}^{b^2}x^{-s}\left\{\int_{\sqrt{x}}^{\min(x/a,b)}f(u)
f\bigl({x\over u}\bigr){\d u\over u}\right\}\d x.
$$

\head
3. The cases  of $\M_k(s)$ when $k=1,3$
\endhead
The analytic continuation of $\M_k(s)$ when $k\le 4$ is
interesting only when $k=1,3$, since $\M_2(s) \equiv \Z_1(s)$,
$\M_4(s) \equiv \Z_2(s)$, and for  $\Z_1(s), \Z_2(s)$ there
is plenty of information (see Section 1). For $k>4$ there
is little information available on $\Z_k(s)$. We have the following

\medskip
THEOREM 4. {\it The function $\M_1(s)$ has analytic continuation to the region
$\s > 0$, where it is regular.
For fixed $\s$ such that ${1\over4} <  \s\le {5\over4}$ it satisfies
$$
\M_1(\s + it) \;\ll_\e\; t^{{3\over4}-\s+\e}(1 + t^{{3\over4}-\s})
\qquad(t \ge t_0 >0).
\leqno(3.1)
$$
We also have, for fixed $\s$ such that $\hf <  \s\le 1$,}
$$
\int_1^T|\M_1(\s + it)|^2\d t \;\ll_\e\; T^{2-2\s+\e},
\leqno(3.2)
$$
$$
\int_1^T|\M_1(\s + it)|^2\d t \;\gg_\e\; T^{2-2\s-\e}.
\leqno(3.3)
$$

\bigskip
THEOREM 5. {\it We have
$$
{\Cal M}_3(s) = \int_1^\infty Z^3(x)x^{-s}\d x = V_1(s) + V_2(s),
$$
say, where $V_2(s)$ is  regular for $\s > 3/4$
and for $\s>1$ the function
$$
V_1(s) = (2\pi)^{1-s}\sqrt{2\over 3}\sum_{n=1}^\infty
d_3(n)n^{-{1\over6}-{2s\over3}}\cos\bigl(3\pi n^{2\over 3}+{\txt{1\over8}}\pi\bigr)
\leqno (3.4)
$$
is regular, where $d_3(n) = \sum_{k\ell m=n}1$.}

\medskip {\bf Proof of Theorem 4.}
To prove the result on the analytic continuation of $\M_1(s)$
we use the author's method of proof [14]. By the use of Laplace
transform of $|\zt|^2$ (see e.g. [28, Theorem 7.15(A)])
it was shown there that
 $\Z_1(s)$ has meromorphic continuation to $\CC$. Thus let
 $$
 {\bar L}(s) := \int_1^\infty Z(y){\roman e}^{-ys}
\d y,\quad L(s) := \int_0^\infty Z(y){\roman e}^{-ys}\ d y\quad(\s = \R s > 0).
$$
Then we have, by absolute convergence, taking initially $\s$ to be sufficiently large
and making the change of variable $xy=t$,
$$\eqalign{&
\int_0^\infty {\bar L}(x) x^{s-1}\d x = \int_0^\infty
\left(\int_1^\infty Z(y){\roman e}^{-xy}\d y\right)
x^{s-1}\d x\cr&
= \int_1^\infty Z(y)\left(\int_0^\infty
x^{s-1}{\roman e}^{-xy}\d x\right)\d y\cr&
= \int_1^\infty Z(y)y^{-s}\d y \int_0^\infty
{\roman e}^{-t}t^{s-1}\d t = \M_1(s)\G(s).\cr}\leqno(3.5)
$$
Since $\G(s)$ has no zeros, it suffices to prove the assertion for
$$\eqalign{&
\int_0^\infty {\bar L}(x) x^{s-1}\d x = \int_0^1{\bar L}(x) x^{s-1}\d x
+ \int_1^\infty{\bar L}(x) x^{s-1}\d x\cr&
= \int_1^\infty{\bar L}(1/x) x^{-1-s}\d x + A(s)\quad\quad(\s>1),\cr}
$$
say, where
$$
A(s) \;:=\; \int_1^\infty{\bar L}_1(x) x^{s-1}\d x
$$
is an entire function. Since
$$
{\bar L}(1/x) =  L(1/x) - \int_0^1 Z(y){\roman e}^{-y/x}
\d y\qquad(x\ge1),
$$
it remains to consider
$$
\eqalign{
\int\limits_1^\infty {\bar L}(1/x)x^{-s-1}\d x &=
\int\limits_1^\infty L(1/x)x^{-s-1}\d x - \int\limits_1^\infty\left(\int_0^1 Z(y)
{\roman e}^{-y/x}\d y\right)
x^{-s-1}\d x \cr&= I_1(s) - I_2(s),\cr}
$$
say. Note that in $I_2(s)$ the integral over $y$ is uniformly bounded, so that
$I_2(s)$ is regular for $\s >0$. To deal with $I_1(s)$ we shall use M. Jutila's
result (see [19, Lemma 2]) that
$$
{\tilde L}(p) \ll 1,\; p = {1\over T} + iu,\; T\ge T_0,
\; 0 \le u \le (T^{1/2}\log T)^{-1},
$$
where
$$
{\tilde L}(p) := \int_0^\infty Z(t)H(\hf + it){\roman e}^{-pt}\d t\qquad(\R p > 0)
$$
with a precisely defined function $H$ which satisfies
$$
H(\hf + it) = 1 + O\left({1\over |t|+1}\right),\quad
H'(\hf + it) =  O\left({1\over (|t|+1)^2}\right).
$$
If we set $k(t) = 1 - H(\hf + it)$, then
$$
I_1(s) = I_3(s) + B(s),
$$
say, where $B(s)$ is regular for $\s>0$ and
$$
\eqalign{
I_3(s) &:= \int_1^\infty\left(\int_0^\infty Z(t)k(t){\roman e}^{-t/x}\d t\right)
x^{-1-s}\d x\cr&
= \int_1^\infty Z(t)k(t)t^{-s}\left(\int_0^\infty{\roman e}^{-u}u^{s-1}\d u\right)\d t
= \G(s)\int_1^\infty Z(t)k(t)t^{-s}\d t.\cr}
$$

Finally note that the author [13] proved that
$$
\Cal I_1(T) \equiv F(T) = \int_1^T Z(y)\d y = O_\e(T^{1/4+\e}),
\leqno(3.6)
$$
which was improved to $F(T) = O(T^{1/4})$ by M. Korolev [22],
who also proved that $F(T) = \Omega_\pm(T^{1/4})$.
M. Jutila [20] gave a different proof of the same
results by establishing precise formulas for $F(T)$.
Integration by parts and (3.6) show that the
$\int_1^\infty Z(t)k(t)t^{-s}\d t$
represents a regular function even for $\s > - 3/4$, implying
that $I_3(s)$, and consequently $\M_1(s)$,
admits analytic continuation to the region $\s >0$, where it is regular.

\medskip
To obtain the pointwise bound (3.1) we use
$$
\M_1(s) = O\Bigl(\frac{1}{t}\Bigr) + \int_{t^{1-\e}}^XZ(x)x^{-s}\d x + \int_X^\infty Z(x)x^{-s}\d x,
\leqno(3.7)
$$
which is valid initially for $\s >1$ and where $X (\gg t)$ is a parameter to be chosen
a little later. One obtains (3.7) by using the argument in (2.12). Integration by
parts and (3.6) show that
$$
\int_X^\infty Z(x)x^{-s}\d x \;\ll_\e\; t^{1+\e}X^{1/4-\s}
\qquad(\s > 1/4,\; X\ll t^C).\leqno(3.8)
$$
The remaining integral in (3.7) is split into $O(\log t)$ integrals of the form
$$
\eqalign{&
\int\limits_Y^{Y'}Z(x)x^{-s}\d x \cr&= 2\int\limits_Y^{Y'}
\sum_{n\le\sqrt{x\over2\pi}}n^{-1/2}
\cos\Bigl(x\log {\sqrt{x/(2\pi)}\over n}- \hf x - \txt{\pi\over8}\Bigr)x^{-s}\d x
+ O\Bigl(\int\limits_Y^{Y'}x^{-1/4-\s}\d x\Bigr),\cr}
$$
where $Y < Y' \le 2Y$, and we used a version of the
classical Riemann--Siegel formula (see e.g., [10, eq.
(4.5)]) for $Z(t)$. Interchanging summation and integration it is seen that
the expression on the right-hand side above is
$$
2\sum_{n\le\sqrt{Y'\over2\pi}}n^{-1/2}\int_{\max(Y,2\pi n^2)}^{Y'}x^{-\s}{\roman
e}^{iF_\pm(x)}\d x + O(Y^{3/4-\s}),\leqno(3.9)
$$
with
$$
\eqalign{
F_\pm(x) &:=  x\log {\sqrt{x/(2\pi)}\over n}- {1\over2} x
- {\pi\over8} \pm t\log x,\cr
F'_\pm(x) &\,= \log {\sqrt{x/(2\pi)}\over n} \pm {t\over x}, \quad
F''_\pm(x) = {1\over2x} \mp {t\over x^2}.\cr}
$$
Consider the contribution of $F_+(x)$, when $F'_+(x) > 0$. If $Y > 4t$ then
$1/(2x) > 2t/(x^2)$, hence by the second derivative test (Lemma 2.1 of [10])
the sum in (3.9) is $\ll Y^{3/4-\s}$. If $Y < t/2$ then $t/(x^2) > 1/x$, hence
again by the second derivative test we obtain a contribution which is
$$
\ll \;Yt^{-1/2}\cdot Y^{1/4-\s} \;\ll\; Y^{3/4-\s}.
$$
If $t/2 \le Y \le 4t$, then $F'_+(x)\gg 1$, hence by the first derivative test
we obtain again a contribution which is $\ll Y^{3/4-\s}$. A
similar analysis  holds for the
contribution of $F_-(x)$, when $F_-''(x) \gg 1/x$. Therefore we have
$$
\int_Y^{Y'}Z(x)x^{-s}\d x \;\ll_\e\; t^{(3/4-\s)(1-\e)} + X^{3/4-\s}.\leqno(3.10)
$$
Choosing $X = t^2$ and noting that $t^{3/2-2\s} \gg t^{-1}$ for $\s \le 5/4$ we
obtain (3.1) from (3.8) and (3.10).

\medskip
{\bf Remark 2}. For $\s>1/2\,$ note that the bound in (3.1) is better than the bound
$$
\Z_1(\s + it) \;\ll_\e\; t^{1-\s+\e}\qquad(0\le \s \le 1,\;t\ge t_0>0),
$$
proved in [18], and for $\s > 2/3$ the bound with the exponent $5/6-\s+\e$ proved
by M. Jutila [19].

\bigskip
The mean square bound (3.2) for $\M_1$ follows by the method of
proof of (see [18, eq. (3.7)])
$$
\int_1^T|\Z_1(\s+it)|^2\d t \;\ll_\e\; T^{2-2\s+\e}\qquad(1/2 \le \s \le 1),
\leqno(3.11)
$$
where instead of Atkinson's formula [1] for the error term in the mean square
formula for $|\zt|$ we use Theorem 1 of M. Jutila [20], which is the analogue
of Atkinson's formula for $Z(t)$, so that there is no need to repeat the
details. In this way it is seen
that for the mean square we do not obtain a better estimate for $\M_1$
than the one derived for $\Z_1$. In fact it was proved (see [12] and [14]) that
$$
\int_1^T|\Z_k(\s+it)|^2\d t \;\gg_\e\; T^{2-2\s-\e}\qquad(k=1,2; \hf < \s \le 1),
\leqno(3.12)
$$
and the lower bound in (3.3) is the analogue of (3.12) for $\M_1$. The proof also
bears similarities to the proofs of (3.12), but we shall give here a sketch of the
proof. From Theorem 2 (with $c = {5\over4},U = X, x \asymp X$) we have
$$
\eqalign{
Z(x) &= {1\over 2\pi i}\int\limits_{{5\over4}-iX}^{{5\over4}+iX}x^{s-1}\M_1(s)\d s
+ O(X^{-1/4})\cr&
= {1\over 2\pi i}\int\limits_{c-iX}^{c+iX}x^{s-1}\M_1(s)\d s +
O\Biggl(\int_c^{5\over4}x^{\s-1}|\M_1(c+iX)|\d\s\Biggr) +  O(X^{-1/4}).\cr}
$$
Now we use the bound (3.1) to obtain that the error terms above are
$$
\ll_\e X^{{1\over2}-c+\e} + X^{\e-{1\over4}} \ll_\e X^{-\e}\qquad(c \ge \hf + 2\e).
$$
Therefore
$$
\int_X^{2X}Z^2(x)\d x \ll \int_X^{2X}\left|
\int_{c-iX}^{c+iX}x^{s-1}\M_1(s)\d s\right|^2\d x + X^{1-2\e}.
$$
Since $Z^2(x) = |\zx|^2$ and $\int_X^{2X}|\zx|^2\d x \gg X\log X$, it follows that
$$
X\log X \ll \int_{X/2}^{5X/2}\f(x)\left|\int_1^Xx^{c+it-1}\M_1(c+it) \d t\right|^2\d x,
\leqno (3.13)
$$
as
$$
\int_0^1 x^{s-1}\M_1(s)\d s \ll 1\qquad (x\asymp X,\; \txt{1\over2} < c \le 1).
$$
Here $\f(x)\,(\ge0)$ is a smooth function supported in $[X/2, 5X/2]$ and equal
to unity in $[X, 2X]$. When we develop the square on
 right-hand side of (3.13) and integrate
sufficiently many times by parts we obtain that
$$
\eqalign{
X\log X& \ll \int_{X/2}^{5X/2}x^{2c-2}\int_1^X\int_{1,|u-t|\le X^\e}^X
|\M_1(c+it)\M_1(c+iu)|\d u\d t\d x\cr&
\ll X^{2c-1}\int_1^X\int_{1,|u-t|\le X^\e}^X \Biggl(|\M_1(c+it)|^2
+ |\M_1(c+iu)|^2\Biggr)\d u \d t\cr&
\ll_\e X^{2c-1+\e}\int_1^X|\M_1(c+it)|^2\d t,\cr}
$$
since the contribution of $|u-t|\ge X^\e$ will be negligibly small.
This implies the assertion (3.3) with $\s = c \ge \hf+2\e$.

\medskip
{\bf Proof of Theorem 5}. Note that
from Theorem 5 of Section 4 (with $k=3$) we obtain (cf. (1.4))
$$
{\Cal I}_3(x) = 2\pi\sqrt{2\over 3}\sum_{n\leqslant ({x\over2\pi})^{3/2}}
d_3(n)n^{-{1\over6}}\cos\bigl(3\pi n^{2\over 3}+{\txt{1\over8}}\pi\bigr)
 + O_\e(x^{3/4+\e}). \leqno(3.14)
$$
Inserting (3.14) in (1.5) we see that
$$
{\Cal Z}_3(s) = V_1(s) + V_2(s),
$$
where $V_2(s)$ (coming from the error term) is obviously regular for $\s > 3/4$
and satisfies $V_2(s) = O(|s|+1)$. Therefore the main problem
is the analytic continuation of
$$
V_1(s) := 2\pi\sqrt{2\over 3}s\int_1^\infty x^{-s-1}
\sum_{n\leqslant ({x\over2\pi})^{3/2}}
d_3(n)n^{-{1\over6}}\cos\bigl(3\pi n^{2\over 3}+
{\txt{1\over8}}\pi\bigr)\d x.\leqno(3.15)
$$
If in (3.15) we invert the order of summation and integration we get
$$\eqalign{
V_1(s) &= -2\pi\sqrt{2\over 3}\sum_{n=1}^\infty
d_3(n)n^{-{1\over6}}\cos\bigl(3\pi n^{2\over 3}+{\txt{1\over8}}\pi\bigr)
\int_{2\pi n^{2/3}}^\infty \d(x^{-s})\cr&
= (2\pi)^{1-s}\sqrt{2\over 3}\sum_{n=1}^\infty
d_3(n)n^{-{1\over6}-{2s\over3}}\cos\bigl(3\pi n^{2\over 3}+{\txt{1\over8}}\pi\bigr).
\cr}\leqno(3.16)
$$
The series in (3.16) converges absolutely for $\s>5/4$.
This is trivial, and we seek a better result. By considering the portion of the
series in (3.16) over $[X,\,2X]$ (for large $X$ and $s = \s+it$ fixed) we want to
show that it is $\ll_\e X^{-\e}$, which provides then the desired analytic continuation
to the right of the $\s$--line. By using the Stieltjes integral representation and then
integration by parts, we are led to two integrals, of which the  relevant one is
$$
J(s;X) := \int_X^{2X}\D_3(x)x^{-{1\over2}-{2s\over3}}
\cos\bigl(3\pi x^{2\over 3}+{\txt{1\over8}}\pi\bigr)\d x.\leqno(3.17)
$$
On applying the truncated Perron inversion formula (see e.g., [10, Appendix]) we have
$$
\D_3(x) = {1\over2\pi i}\int_{{1\over2}-iX}^{{1\over2}+iX}{\z^3(w)x^w\over w}\d w
+ O_\e(X^\e)\qquad(X \le x \le 2X),\leqno(3.18)
$$
where as usual $\D_3(x)$ is the error term in the asymptotic formula for
the summatory function of $d_3(n)$.
The error term in (3.18) contributes to the integral in (3.17)
$\;\ll_\e X^{{1\over2}-{2\s\over3}+\e}
\ll_\e X^{-\e}$ for $\s > 3/4$. The main term in (3.18) produces
$$
{1\over2\pi i}\int_{{1\over2}-iX}^{{1\over2}+iX}{\z^3(w)\over w}
\left(\int_X^{2X}x^{-2\s/3}\exp(iF_\pm(x))\d x\right)\d w,
$$
where
$$
w = \hf + iv, \;s = \s + it,\; F_\pm(x) :=
\bigl(v - (2t)/3\bigr)\log x \pm 3\pi x^{2/3}.
$$
Note that the saddle point (root of the equation $F_\pm'(x) = 0$)
$$
\;x_0 = \left({|v-(2t)/3|\over2\pi}\right)^{3/2}\in [X,2X]
\qquad({\roman {for}} \;v \asymp X^{2/3}),
$$
in which case  $|F_\pm''(x_0)|^{-1/2} \asymp X^{2/3}$. Hence
by the saddle-point method the total contribution to (3.18) is
$\ll_\e X^{(2/3)(1-\s)+\e}$,
and this provides the desired analytic continuation of ${\Cal Z}_3(s)$ only to
$\s >1$ as before. One can make the calculation of (3.17) simpler by making
the change of variable $x^{2/3}=y$, after $\D_3(x)$ is replaced by (3.18).
However at present I do not see any better way to tackle the
problem of the  analytic
continuation of ${\Cal Z}_3(s)$, although I feel that it can be done.

\medskip
{\bf Remark 3}. It is curious that obviously the shapes of ${\Cal M}_k(s)$ for $k = 1,2,3,4$
(the cases when we know something relevant) are totally different! The fact that $Z(x)$ is an
oscillating function, while $|\zx|$ is non-negative is reflected in what we expect:
${\Cal M}_{2\ell}(s) ={\Cal Z}_\ell(s) $ should have a pole of order $\ell^2+1$ at
$s=1$, while ${\Cal M}_{2\ell-1}(s)$ should be regular at $s=1$, at least for $1 \le \ell
\le 4$.

\medskip

\head
4. The cubic moment of $Z(t)$
\endhead
Let, as usual, $d_k(n)$ denote the number of ways $n$ can be written as
a product of $k$ factors, so that $d_k(n)$ is the multiplicative function
generated  by $\z^k(s)$. In particular, $d_1(n)\equiv 1$ and
$d_2(n) \equiv d(n)$, the number of divisors of $n$.
  To prove the second part of Theorem 4 we need the case $k=3$ of

\medskip
 THEOREM 6. {\it For fixed $k = 1,2,3,4$ we have
$$
\eqalign{
\int\limits_T^{2T}Z^k(t)\d t&= 2\pi\sqrt{2\over k}
\sum_{({T\over 2\pi})^{k/2}\le n\le ({T\over\pi})^{k/2}}
d_k(n)n^{-{1\over2}+{1\over k}}\cos\bigl(k\pi n^{2\over k}+{\txt{1\over8}}(k-2)\pi\bigr)
\cr&
+ \ldots +O_\e(T^{k/4+\e}),\cr}\leqno(4.1)
$$
where $+\ldots$ denotes terms similar to the one on the right-hand side of }(4.1),
{\it with the similar cosine term, but of a lower order of magnitude}.

\medskip
{\bf Proof of Theorem 6}.
For $Z^k(t)$ we shall use a finite, smoothed sum,
which is a form of the so-called approximate functional equation.
One could also use
a form  of the approximate functional equation which comes from the
so-called ``reflection principle'' (see e.g., Chapter 4 of [10]).
However, to have a symmetric expression we shall use essentially
a variant  of the approximate functional equation for $\z^k(s)$ which
is to be found in Chapter 4 of [11]. To this end let $\rho(x)$
be a non-negative, smooth function
supported in $\,[0,2]\,$, such that $\rho(x) = 1$ for $0 \le x\le 1/b$
for a fixed constant $b>1$, and $\rho(x) + \rho(1/x) = 1$ for all $x$
(an explicit construction of $\rho(x)$ was given in [11]).
Let $\tau = \tau(k,t)$ be defined as
$$
\log\tau = - k{\chi'(\hf+it)\over\chi(\hf+it)}.\leqno(4.2)
$$
We write
$$
\chi(s) = \pi^{s-1/2}{\Gamma(\hf - \hf s)\over\Gamma(\hf s)}
= \left({2\pi\over t}\right)^{\s+it-1/2}{\roman e}^{i(t+\pi/4)}
\left(1+ O\Bigl({1\over t}\Bigr)\right)\leqno(4.3)
$$
by using Stirling's formula for the gamma-function. Here
$s = \s+it, 0\leqslant \s\leqslant 1, t \geqslant t_0 > 0,$
and note that the $O$-term in (4.3) admits an asymptotic expansion
in terms of negative powers of $t$. Therefore
$$
{\chi'(\hf+it)\over\chi(\hf+it)} = \log\left({2\pi\over t}\right)
+ O\left({1\over t^2}\right),
$$
and we obtain
$$
\tau = \left({t\over2\pi}\right)^k\left(1 + O\left({1\over t^2}\right)\right),
\leqno(4.4)
$$
and again the $O$-term in (4.4) admits an asymptotic expansion
in terms of negative powers of $t$. In the course of the proof of Theorem 5.2 of [11]
it was shown that ($1\ll x,y \ll t^k, xy = \tau,
s =\s+it, t \geqslant t_0>0, 0\leqslant \s \leqslant1$)
$$
\z^k(s) = \sum_{n=1}^\infty d_k(n)\rho\Bigl({n\over x}\Bigr)n^{-s}+\chi^k(s)
\sum_{n=1}^\infty d_k(n)\rho\Bigl({n\over y}\Bigr)n^{s-1} + R_k(t),
\leqno(4.5)
$$
say, where for any fixed $A>0$
$$
R_k(t) \;\ll_\e\; t^{-A} + t^{\e-1}\int_{-t^\e}^{t^\e}|\z(\s+it - \e+iv)|^k\d v.
\leqno(4.6)
$$
Thus from (4.3)--(4.6) we obtain, with $b=2, \s=\hf, x=y=\sqrt{\tau},
t\geqslant t_0>0$ the following

\medskip
LEMMA 2. {\it We have
$$
\eqalign{
\int_T^{2T}Z^k(t)\d t&= 2\int_T^{2T}\sum_{n\leqslant2\sqrt{\tau}}
\rho\Bigl({n\over \sqrt{\tau}}\Bigr)d_k(n)n^{-1/2}\cos F_k(t)\d t
+\ldots\cr&
+ O\left(T^{\e-1}\int_{T/2}^{5T/2}|\zt|^k\d t\right),\cr}\leqno(4.7)
$$
where $\tau$ is given by }(4.2) {\it and} (4.4),
$+\ldots$ {\it denotes terms similar to the one on the right-hand side of }(4.7),
{\it but of a lower order of magnitude, and where
}
$$
F_k(t) := t\log\Biggl\{{\Bigl({t\over2\pi}\bigr)^{k/2}\over n}\Biggr\} -
{kt\over 2} - {k\pi\over8}. \leqno(4.8)
$$

\medskip
To evaluate the left-hand side of (4.7) we write first
$$
\eqalign{&
2\int_T^{2T}\sum_{n\leqslant2\sqrt{\tau}}
\rho\Bigl({n\over \sqrt{\tau}}\Bigr)d_k(n)n^{-1/2}\cos F_k(t)\d t \cr&
= 2\sum_{n\le T_0}d_k(n)n^{-1/2}\R\Bigl\{
\int_{T_1}^{2T}\rho\Bigl({n\over \sqrt{\tau}}\Bigr){\roman e}^{iF_k(t)}\d t\Bigr\}.\cr}
\leqno(4.9)
$$
Here
$$
T_0 = 2\sqrt{\tau(k,2T)} = 2\left({T\over\pi}\right)^{k/2}
\left(1 + O\Bigl({1\over T^{2}}\Bigr)\right),
\qquad T_1 = \max\Bigl(T, \tau^{-1}(k,(n/2)^2\Bigr),
$$
where $\tau^{-1}$ is the inverse function of $\tau$, so that
$$
\tau^{-1}(k,(n/2)^2) = 2\pi\left({n\over2}\right)^{2/k}
\left(1 + O\Bigl({1\over T^{2}}\Bigr)\right).
$$
Now we split the range of summation over $n$ on the right hand side of (4.9) as follows.
Let
$$
\eqalign{
I_1 &:= \left[1, \,\left({T\over2\pi}\right)^{k/2} - T^{k/2-1/2+\e}\right),\cr
I_2 &:= \left[\left({T\over2\pi}\right)^{k/2} - T^{k/2-1/2+\e},
\left({T\over2\pi}\right)^{k/2} + T^{k/2-1/2+\e}\right),\cr
I_3 &:= \left[\left({T\over2\pi}\right)^{k/2} + T^{k/2-1/2+\e},
\left({T\over\pi}\right)^{k/2} - T^{k/2-1/2+\e}\right],\cr
I_4 &:= \left(\left({T\over\pi}\right)^{k/2} - T^{k/2-1/2+\e},
\left({T\over\pi}\right)^{k/2} + T^{k/2-1/2+\e}\right],\cr
I_5 &:= \left(\left({T\over\pi}\right)^{k/2} + T^{k/2-1/2+\e}, T_0\right].\cr}\leqno(4.10)
$$
In the integrals over where $n\in I_1$ and $n\in I_5$ we integrate by parts, writing
$$
\int \rho\Bigl({n\over \sqrt{\tau}}\Bigr){\roman e}^{iF_k(t)}\d t
= \int{\rho\Bigl({n\over \sqrt{\tau}}\Bigr)\over i\log\bigl\{(t/2\pi)^{k/2}/n\bigr\}}
\d {\roman e}^{iF_k(t)}.\leqno(4.11)
$$
Note that the derivatives of $\rho\Bigl(n/\sqrt{\tau}\Bigr)$, considered as a function of $t$,
decrease after each integration by parts  by a factor of $t$,
while in $\sum_{n\in I_1}\int$ we have
$$\eqalign{
{\left(1\over \log\bigl\{(t/2\pi)^{k/2}/n\bigr\}\right)}^{'}
&= - {2\over kt\log^2\bigl\{(t/2\pi)^{k/2}/n\bigr\}} \cr&\ll_\e
{1\over T\log^2\left\{CT^{k/2}\over T^{k/2} +O(T^{k/2-1/2+\e})\right\}}
\ll_\e T^{-2\e}.\cr}\leqno(4.12)
$$
Therefore if we integrate by parts sufficiently many times, the contribution
will be negligible. The sums over the integrated terms are essentially partial
sums of $\z^k(\hf+iu), u\asymp T$, when we remove the monotonic coefficients
$\rho$ from the sums over $n$ by partial summation.
The resulting sums are bounded by Perron's inversion formula
(see e.g., the Appendix of [10]).  Since $\zt \ll t^c$
for some $c< 1/6$ (ibid., Chapter 7), we see that
$$
\sum\limits_{n\in I_1}\; +\; \sum\limits_{n\in I_5} \;\ll\; T^{k/6}.\leqno(4.13)
$$
Note that (cf. (4.8))
$$
F_k'(t) = \log\bigl\{(t/2\pi)^{k/2}/n\bigr\},\quad F_k''(t) = k/(2t).\leqno(4.14)
$$
The integrals when $n\in I_2 \cup I_4$ are estimated as $\ll T^{1/2}$
by the second derivative test (see Chapter 2 of [10]), and then trivial
estimation gives
$$
\sum\limits_{n\in I_2}\; +\; \sum\limits_{n\in I_4} \;\ll_\e\;
T^{1/2}T^{k/2-1/2+\e}T^{-k/4} = T^{k/4+\e}.\leqno(4.15)
$$
Finally when, in (4.9), we have $n\in I_3$, then the saddle point (root of $F_k'(t) = 0$),
namely
$$
t_0 \;\equiv\; c_n\; :=\; 2\pi n^{2/k}\leqno(4.16)
$$
lies in $[T_1,\, 2T]$.
For $\int_{T_1}^{2T}$ we could use a general result on exponential integrals,
such as the following [21, Lemma III.2], which says that
$$
\eqalign{&
\int_a^b \f(x)\exp\Bigl(2\pi if(x)\Bigr)\d x = {\f(c)\over\sqrt{f''(c)}}
{\roman e}^{2\pi if(c)+\pi i/4} + O(HAU^{-1})\cr&
+ O\bigl(H\min(|f'(a)|^{-1},\sqrt{A}\,\bigr) +
O\bigl(H\min(|f'(b)|^{-1},\sqrt{A}\,\bigr)\cr}
\leqno(4.17)
$$
if $f'(c) = 0$, $a \le c \le b$, and the following conditions hold: $f(x) \in
C^4[a,b]$, $\f(x) \in C^2[a,b]$, $f''(x) > 0$ in $[a,b]$, $f''(x) \asymp
A^{-1}$, $f^{(3)}(x) \ll A^{-1}U^{-1}$, $f^{(4)}(x) \ll A^{-1}U^{-2}$,
$\f^{(r)}(x) \ll HU^{-r}\;(r = 0,1,2)$ in $[a,b],\;0 < H,A < U, \,
0 < b - a \le U$. In our case $f(x) = {1\over2\pi}F_k(x)$,
so that $f''(c) = k/(4\pi c)$, and

\smallskip

$$
{\f(c_n)\over\sqrt{f''(c_n)}}
{\roman e}^{2\pi if(c_n)+\pi i/4} = \pi \sqrt{2\over k}\,
n^{1\over k}\exp\left(-k\pi in^{2\over k}
+{(2-k)\pi i\over8}\right)\Bigl\{1+O\bigl({1\over T^2}\bigr)\Bigr\}.\leqno(4.18)
$$
But, as remarked in [13], in our case the last two
error terms in (4.17) are large, and thus
it is more expedient to carry out the evaluation by the saddle point technique
directly, that is, by using a suitable contour in the complex plane.

To this end, if $T_1 = T$ (the other case is similar) let
$\Cal L_1$ be the segment $T
- iv \,(0 \le v \le {1\over\sqrt{2}}T^{1-\e})$,
$\Cal L_2$ is the segment $x - i{1\over\sqrt{2}}T^{1-\e}\,(0\le x
\le c_n - {1\over\sqrt{2}}T^{1-\e})$,  $\Cal L_3$ is the segment
$ c_n + v{\roman e}^{{1\over4}\pi i}$,
$-{1\over\sqrt{2}}T^{1-\e} \le v \le {1\over\sqrt{2}}T^{1-\e}$,
$\Cal L_4$ is the segment $x +  i{1\over\sqrt{2}}T^{1-\e}\,
(c_n + {1\over\sqrt{2}}T^{1-\e}\le x\le 2T)$, and finally $\Cal L_4$
is the segment joining the points $2T + i{1\over\sqrt{2}}T^{1-\e}$
and $2T$.

As a simplification we develop $\rho(n/\sqrt{\tau})$ by Taylor's formula at the point
$t_0 = c_n$ when $t\in [c_n-T^{1-\e}, c_n+T^{1-\e}]$, and at other appropriate points for
other values of $t$. An alternative approach is to use the Mellin inversion formula:
$$
\rho(x) = {1\over2\pi}\int_{d-i\infty}^{d+i\infty}R(s)x^{-s}\d s\;\;(d>0),
\quad R(s) = \int_0^\infty \rho(x)x^{s-1}\d x.
$$
The function $R(s)$ is odd, and of fast decay.

As already noted the derivatives
of $\rho(n/\sqrt{\tau})$,
considered as a function of $t$,
decrease each time by a factor of $t$. Since the length of the interval
is $2T^{1-\e}$, it is possible to take finitely many terms in Taylor's
formula so that the total contribution of the error term is negligible,
namely $\ll 1$. The remaining integrals will be all of the
same type, with the same exponential factor, and the largest one will
be the first one, namely the one with ($c_n = 2\pi n^{2/k}$)
$$
\rho\left({n\over\sqrt{\tau(k,c_n)}}\right) = \rho\left({n\over n(1+O(T^{-2})}\right)
= \rho(1) + O(T^{-2}) = \frac{1}{2} + O(T^{-2}),
$$
since $\rho(x) + \rho(1/x)=1$. Here actually the $O$-term above has an asymptotic expansion.
After that
we replace the subinterval integral over  $[T_1, 2T]$,
by Cauchy's theorem, by $\cup_{j=1}^5\int_{{\Cal L}_j}$. Therefore

$$\eqalign{&
2\sum_{n\in I_3} d_k(n)n^{-1/2}\R\Bigl\{\int_{T_1}^{2T}
\rho\Bigl({n\over \sqrt{\tau}}\Bigr){\roman e}^{iF_k(t)}\d t\Bigr\}\cr&
= 2\sum_{n\in I_3} d_k(n)n^{-1/2}\R\Bigl\{\bigcup_{j=1}^5\int_{{\Cal L}_j}
{\roman e}^{iF_k(z)}\d z\Bigr\} + \ldots\,,\cr}\leqno(4.19)
$$
where $+ \ldots$ has the same meaning as before.
On $\Cal L_3$ we have (since $F_k'(c_n) = 0$)
$$
iF_k(z) = iF_k(c_n) + i{v^2\over2!}{\roman e}^{{1\over2}\pi i}F_k''(c_n)
+ i{v^3\over3!}{\roman e}^{{3\over4}\pi i}F_k'''(c_n)
+ i{v^4\over4!}F_k^{(4)}(c_n) + \cdots\,. \leqno(4.20)
$$
Note that, since $v\ll T^{1-\e}$,
$$
v^{m}F_k^{(m)}(c_n) \;\ll_{m,\e}\; T^{m(1-\e)}T^{1-m} \;=\; T^{1-m\e}\qquad(m>1).
\leqno(4.21)
$$
Hence if we choose $M = M(k,\e)$ sufficiently large, then (4.7) shows that the terms
of the series in (4.20) for $m > M$,
on using $\exp z = 1 + O(|z|)$ for $|z| \le 1$,
will make a negligible contribution. We have
$$
\exp(iF_k(z)) = \exp(iF_k(c_n))\exp(-\hf v^2 F''(c_n))
\exp\left(\sum_{m=3}^M d_m v^m F^{(m)}(c_n)\right)
$$
with  $d_m = \exp((m+2){\pi i\over4})/m!$.
The last exponential factor is expanded by Taylor's series, and again
the terms of the series (with $v^m$) for large $m$ will make a negligible
contribution. In the remaining terms we restore integration over
$v$ to the whole real line, making a very small error. Then we use
the classical integral (see e.g., the Appendix of [10])
$$
\int_{-\infty}^\infty \exp(Ax - Bx^2)\,\d x \;=\;
\sqrt{\pi\over B}\exp\left({A^2\over4B}\right)\qquad(\R B > 0).\leqno(4.22)
$$
By differentiating (4.22) as a function of $A$ we may explicitly
evaluate integrals of the type
$$
\int_{-\infty}^\infty x^{2m}\exp(-Bx^2)\,\d x\qquad(\R B > 0,\;
m = 0,1,2,\ldots\,).
$$
It transpires that the largest contribution ($= \sqrt{\pi}$) will come
from the integral with $m = 0$, which will coincide with the contribution
of the main term in (4.17).

It remains to deal with the remaining integrals over $\Cal L_j$.
The integrals over $\Cal L_1$ and $\Cal L_5$, and likewise the
integrals over $\Cal L_2$ and $\Cal L_4$ are estimated analogously.
On $\Cal L_4$ we have
$$
z = x+iH, \, c_n + {H\over\sqrt{2}} \le x \le 2T\,,\, H = T^{1-\e}.
$$
On using Taylor's formula we obtain
$$
\exp(iF_k(z)) = \exp(iF_k(x) - i{H^2\over2!}F_k''(x) +\ldots)
\exp(-HF_k'(x) + {H^3\over3!}F_k'''(x) -\ldots).
$$
Similarly as in (4.21) it follows that we may truncate the series after a finite
number of terms with a negligible error. Observe that the real-valued term in the
exponential is negative, and that the derivative of the imaginary part is dominated
by
$$\eqalign{
F_k'(x) &= \log{\bigl({x\over2\pi}\bigr)^{k/2}\over n }\ge
\log{\Bigl(n^{2/k} + H/\sqrt{2}\Bigr)^{k/2}\over n}\cr
&= \log\Bigl(1 + {H\over\sqrt{8}\pi n^{2/k}}\Bigr)^{k/2} \ge A_kHT^{-1} = A_kT^{-\e}\cr}
$$
for some constant $A_k>0$. Hence by the first derivative test the total contribution
of such terms is
$$
\ll_\e\; T^{k/4+\e}.\leqno(4.23)
$$

On $\Cal L_5$ we have $z = 2T + iy, 0\le y\le H, H=T^{1-\e}.$ This gives
$$
iF_k(z) = iF_k(2T) - yF_k'(2T) - i{y^2\over2!}F_k''(2T) + {y^3\over 3!}F_k'''(T) - \ldots\,,
$$
where, as before, we may truncate the series after a finite
number of terms with a negligible error. Therefore the integral over $\Cal L_5$ becomes
$$
i{\roman e}^{iF_k(2T)}\int_0^H {\roman e}^{f(y)}{\roman e}^{ig(y)}\d y,
$$
say, with real-valued
$$
f(y) = - yF_k'(2T)  + {y^3\over 3!}F_k'''(T)\ldots\,,\;\;
g(y) = - {y^2\over2!}F_k''(2T) + {y^4\over4!}F_k^{(4)}(2T) +\ldots\,.
$$
Then we have
$$
\int_0^H = \int_0^{\sqrt{T}} + \int_{\sqrt{T}}^H = J_1 + J_2,
$$
say. We write $J_1$ as
$$
J_1 = -{1\over F_k'(2T)  + {y^2\over 2!}F_k'''(2T)\ldots\,}
\int_0^{\sqrt{T}}{\roman e}^{ig(y)}\d\Bigl({\roman e}^{f(y)}\Bigr)
$$
and integrate by parts. We obtain the same type of exponential integral, only
smaller by a factor of
$$
\ll y{F_k''(2T)\over F_k'(2T)} \ll T^{1/2}\cdot{1\over T}\cdot T^{1/2-\e} = T^{-\e},
$$
since
$$
 F_k'(2T) \ge \log {(T/\pi)^{k/2}\over n} \ge \log {(T/\pi)^{k/2}\over
 (T/\pi)^{k/2}- T^{k/2-1/2+\e}} \ge T^{\e-1/2}.
$$
This means that, after sufficiently many integrations by parts, the ensuing integral
will be negligible, while the integrated terms will  be $\ll T^{k/4+\e}$ as in (4.9).
Finally in $J_2$
$$
yF_k'(2T)  - {y^3\over 3!}F_k'''(T)\ldots \;\ge\; CT^{1/2}\cdot T^{\e-1/2} = CT^{\e},
$$
so that  ${\roman e}^{f(y)}$ is negligibly small. The net result of our considerations
is that in the evaluation of the right-hand side of (3.1)
the main terms, arising from the saddle point terms, are given by (4.18), while
all the error terms are $\ll_\e T^{k/4+\e}$.

\medskip
{\bf Remark 4}. With a more careful analysis one can get
rid of the terms implied by $+\ldots$ in (4.1). The same also follows
if one uses an idea of Prof. Matti Jutila, who
kindly informed me that the above proof may
be simplified as follows. The method may be traced back to
E.C. Titchmarsh [28, p. 261], and a sketch is as follows. Note that
$$
\eqalign{&
\int_T^{2T}Z^k(t)\d t = -i\int_{{1\over2}+iT}^{{1\over2}+2iT}\chi^{-k/2}(s)\z^k(s)\d s
\cr&
= -i\left(\int_{1+\e+iT}^{1+\e+2iT} + \int_{{1\over2}+iT}^{1+\e+iT}
-\int_{{1\over2}+\e+iT}^{1+\e+2iT}\right)\chi^{-k/2}(s)\z^k(s)\d s.\cr}
$$
On $\s = 1+\e$ we have $\z^k(s) = \sum_{n=1}^\infty d_k(n)n^{-s}$, so that
the above expression is seen to be
$$
\sum_{n=1}^\infty d_k(n)n^{-1-\e}\int_T^{2T}\Bigl({t\over2\pi}\Bigr)^{{k\over4}
+ {k\e\over2}}{\roman e}^{iF_k(t)}\d t + O_{\e,k}(T^{k/4+\e})
$$
for $k\le 4$. The exponential integral is evaluated by (4.17), and Theorem 6
will follow. I am grateful to Prof. Jutila for pointing this out to me.

\head
5. Natural boundaries
\endhead

If a Dirichlet series $F(s)$ has a (meromorphic) continuation to
$\R s > \s_0$, then the line $\R s = \s_0$ is said to be the
{\it natural boundary} of $F(s)$ if the poles of $F(s)$ are dense
on $\R s = \s_0$, so that $F(s)$ cannot be continued analytically
to  $\R s \le \s_0$.
The history of natural boundaries for Dirichlet series goes at
least back to T. Estermann [9]. For example, one has
$$
\eqalign{
\sum_{n=1}^\infty d_k^2(n)n^{-s} &= \zeta^{k^2}(s)\prod_pP_k(p^{-s})
\qquad(\R s >1),\cr
P_k(u) &:= (1-u)^{2k-1}\sum_{n=0}^k{k-1\choose n}^2u^k,\cr}
$$
and Estermann showed that the  above Euler product has meromorphic
continuation to $\R s>0$, but has the line $\R s =0$ as the natural
boundary when $k>2$. In fact, his result holds for a class of Dirichlet series
of which the above product is a special case. Estermann's results were
generalized by G. Dahlquist [5], and recent investigations include the works
of G. Bhowmik and J.-C. Schlage--Puchta [2], [3].

\smallskip
If $A(s) = \sum_{n=1}^\infty a_nn^{-s}$ in its region of absolute convergence
$\R s > \s_a$, then by Perron's inversion formula
$$
\sum_{n\le x}a_n = {1\over2\pi i}\int_{c-i\infty}^{c+i\infty}
A(s){x^s\over s}\d s\qquad(x\not \in\NN,\; c> \s_a).\leqno(5.1)
$$
In practice one wants to shift the line of integration in (5.1) to the left,
to reduce the contribution of the term $x^s$. This is possible only if $A(s)$
is holomorphic on the new path. If $\s = \s_0\;(< \s_a)\,$ is the natural boundary of
$F(s)$, then we cannot have $c\le \s_0$, hence the usefulness of (5.1) is
limited if $\s_0$ exists. This is one of the reasons which makes the study
of natural boundaries of Dirichlet series important.

\smallskip
The interest in natural boundaries for $\Z_k(s)$ begins with the notes of
A. Diaconu [7], followed by the author's notes [16], and the papers of Diaconu,
Garrett, Goldfeld  [8] and Y. Motohashi [25]. Note that $\Z_k(s)$
does not have an Euler product, which makes the problem more difficult.
It is conjectured in all these works  that  the analytic continuation of
$\Z_3(s)\,(\equiv\,\M_6(s))$ has $\R s = \hf$ as the natural boundary,
and that, more generally, $\Z_k(s)$ for $k\ge3$ has
$\R s = \hf$ as the natural boundary.
A full proof of this important claim concerning $\Z_k(s)$ would be most welcome.
The basic idea that leads to it is
simple, and is open to  generalizations. Namely on p. 2 of [6]
(or p. 3 of [7]) it is said
that the analytic continuation of ($s,w$ are complex variables)
$$
\int_1^\infty {\Bigl({m\over n}\Bigr)}^{ix}L(s_1+ix)L(s_2-ix)\,x^{-w}\d x\leqno(5.2)
$$
produces the analytic continuation of
$$
\int_1^\infty |F(\s +ix)L(\s+ix)|^2\,x^{-w}\d x, \quad F(s) =
\sum_{n=1}^\infty f(n)n^{-s}\leqno(5.3)
$$
under some reasonable conditions, simply by squaring out $|F|^2$ and summing
over the relevant $m,n$. In Proposition 2.6 on p. 3 this approach is discussed
when $L$ is the zeta-function of a holomorphic cusp form of weight $\k$
for $SL(2,\,\ZZ)$.

If in (5.3) we take $F = \z, L = \z^2, \s = \hf$, then we have to observe
that $\Z_2(s)$ has (see Section 1) infinitely many poles
at $s = {1\over2} \pm i\k_j\,\left(\k_j =\sqrt{\lambda_j -
{1\over4}}\,\right)$.
Heuristically, when we sum over various $m,n$ in (5.2) to get the analytic
continuation of $\Z_3(w)$, each of the poles ${1\over2} \pm i\k_j$ will
be somewhat perturbed. Their totality will be dense on the $\hf$--line,
and will produce the $\hf$--line as the natural boundary for $\Z_3(w)$.
Inasmuch as this seems plausible, a rigorous proof is in order.

\smallskip
Suppose that one has found the analytic continuation of $\Z_3(w)$ to the
right of the $\hf$--line. Then it is seems plausible that $\Z_3(w)$
(being more complex that $\Z_2(w)$) will have infinitely many poles
as well. Where are these poles located? One does not expect
them be too near the $\hf$--line,
so the ${3\over4}$--line is a very good candidate to contain
infinitely many poles of $\Z_3(w)$. But by
the principle inherent in (5.2)--(5.3), then the ${3\over4}$--line would be
a natural barrier for $\Z_4(w)$, and so on -- each $\Z_k(w)$ would,
with increasing $k$, have poles nearing the $1$--line.

\smallskip
The recent work of Y. Motohashi [25] on (5.3) (when $L = \z^2$)
supports the claim that $\Z_3(s)$ has $\s = 1/2$ as the natural boundary.
The author says: ``Our
theorem suggests that the Mellin transform $\int_1^\infty |\zx|^6x^{-s}\d x$ should
have the line $\R s = 1/2$ as a natural boundary... The same was also speculated also
by a few people other than us, but it appears that our theorem is so far the sole
explicit evidence supporting the observation.''

The natural boundary of $\Z_3(s)$ on $\R s = 1/2$ indicates certainly
a complicated structure of the error term $E_{3}(T)$
for the sixth moment of $|\zt|$, but in itself it does not exclude
the possibility of the bound $E_{3}(T) \ll_\e T^{1/2 + \e}$.
If $E_{3}(T) \ll_\e T^{\theta + \e}$ with $\theta$
as small as possible, then $\Z_3(s)$ would have singularities on
$\R s = \theta$, if $1/2< \theta <1$. Inasmuch as it seems plausible (to me) that
$\theta = 3/4$, this is a major unsolved problem.

\smallskip
{\bf Remark 5}. I believe that ($P_9(y)$ is an explicit polynomial of degree nine)
$$
\eqalign{&
\int_0^T|\zt|^6\d t = TP_9(\log T) + E_3(T),\cr& E_3(T)
= O_\e(T^{3/4+\e}),\quad E_3(T) = \Omega(T^{3/4})\cr}
$$
holds, where the main term $TP_9(\log T)$ is the one predicted by Conrey et al. [4].
However in [4] the error term is indicated to be (in all cases) $O_\e(T^{1/2+\e})$,
which I do not think can be true.

\smallskip
In what concerns the true order of higher moments of $|\zt|$, the situation
is even more unclear.
Already for the eighth moment it is hard to ascertain what goes on, much less for the
higher moments. The main term for the general
$2k$-th moment should involve a main
term of the type suggested by [2], but it could turn out that the error term
$$
E_k(T)  = \int_0^T|\zt|^{2k}\d t - TP_{k^2}(\log T) \qquad(k\in \NN)
$$
in the general case (when $k \ge 4$) contains expressions which make
it {\it larger} than the   term $TP_{k^2}(\log T)$. For this see the discussion in [12]
(also [24, pp. 218-219]). Essentially the argument is as follows.
In general, from the knowledge about the order of $E_k(T)$
one can deduce a bound for $\zeta({1\over 2} + iT)$ via the estimate
$$
\zeta(\hf + iT) \;\ll\; (\log T)^{(k^2+1)/(2k)} +
{\Bigl(\,\log T \max_{t\in [T-1,T+1]} \vert E_k(t)\vert\,\Bigr)}^{1/(2k)},
\leqno(5.5)
$$
which is Lemma 4.2 of [11]. The conjectured bounds
$$
E_k(T) \;\ll_\e\; T^{k/4+\e}  \qquad(k\le4)\leqno(5.6)
$$
all imply $\zt \ll_\e |t|^{1/8+\e}$, which is out of reach at present, but
is still much weaker than the Lindel\"of hypothesis that
$\zt \ll_\e |t|^{\e}$. On the other hand, we know that the omega-result
$$
E_k(T) \;=\; \Omega(T^{k/4})  \leqno(5.7)
$$
hold for $k =1,2$, and as already explained, there are reasons to believe
that (5.7) holds for $k=3$. Perhaps it holds for $k=4$ also, but the
truth of (5.7) for any $k>4$ would imply that the Lindel\"of hypothesis
is false, and {\it ipse facto} the falsity of the Riemann hypothesis
(that all complex zeros of $\z(s)$ satisfy $\R s = 1/2$). Namely
it is well-known (see e.g., [10] or [28]) that the Riemann hypothesis
implies even $\log|\zt|\ll \log |t|/\log\log |t|$, which is stronger than
the Lindel\"of hypothesis ($\Leftrightarrow \log|\zt|\ll_\e \e\log |t|$).
The reason why, in general, (5.7) makes sense
is that a bound $E_k(T) \ll T^{c_k}$ for some fixed $k\;(>4)$ with $c_k < k/4$ would
imply (by (5.5)) the bound  $\zt \ll_\e |t|^{c_k/(2k)+\e}$ with $c_k/(2k) < 1/8$.
But the most one can get (by using (5.5)) from the error term in the mean square
and the fourth moment of $|\zt|$ is the bound
$$
\zt \;\ll_\e\; |t|^{1/8+\e}.
$$
It does not appear likely to me that, say from
the twelfth moment ($k=6$), one will get a better
pointwise estimate for $\zt$ than what one
can get from the mean square formula ($k=1$).
Nothing, of course, precludes yet that this does not happen,
just that it appears to me not to be likely. As in all such dilemmas,
only rigorous proofs will reveal in due time the real truth.

\smallskip

\Refs

\item{[1]} F.V. Atkinson,  The mean value of the Riemann zeta-function,
 Acta Math. {\bf81} (1949), 353-376.

\item{[2]} G. Bhowmik and J.-C. Schlage-Puchta, Natural boundaries of Dirichlet
series, Func. Approx. Comment. Math. {\bf37}(2007), 17-29.

\item{[3]} G. Bhowmik and J.-C. Schlage-Puchta, Essential singularities of Euler
products, to appear, see {\tt arXiv:1001.1891}.

\item{[4]} J.B. Conrey, D.W. Farmer, J.P. Keating, M.O. Rubinstein
and N.C. Snaith, Integral moments of $L$-functions,
Proc. London Math. Soc. (3) {\bf91}(2005), 33-104.

\item{[5]} G. Dahlquist, On the analytic continuation of Eulerian products,
Ark Mat. {\bf1}(1952), 533-554.

\item{[6]} A. Diaconu, The function $\Z_3(w)$ has natural boundary,
Notes of October 24, 2006.

\item{[7]} A. Diaconu, P. Garrett and D. Goldfeld, Natural
boundaries and a correct notion of integral moments of $L$-functions,
preprint, 2009.

\item{[8]} A. Diaconu, D. Goldfeld and J. Hoffstein,
Multiple Dirichlet series and moments of zeta and $L$-functions
Compos. Math. {\bf139}(2003), 297-360.

\item{[9]} T. Estermann, On certain functions represented by
Dirichlet series, Proc. London Math. Soc. {\bf27}(1926), 435-448.

\item{[10]} A. Ivi\'c, The Riemann zeta-function, John Wiley \&
Sons, New York 1985 (2nd edition. Dover, Mineola, New York, 2003).

\item{[11]} A. Ivi\'c,  The mean values of the Riemann zeta-function,
LNs {\bf 82}, Tata Inst. of Fundamental Research, Bombay 1991 (distr. by
Springer Verlag, Berlin etc.).

\item {[12]} A. Ivi\'c, On some conjectures and results for the Riemann zeta-function
and Hecke series, Acta Arithmetica  {\bf109}(2001), 115-145.

\item {[13]} A. Ivi\'c, On the integral of Hardy's function, Arch. Mathematik
{\bf83}(2004), 41-47.

\item {[14]} A. Ivi\'c, The Mellin transform of the square of Riemann's zeta-function,
International J. of Number Theory {\bf1}(2005), 65-73.

\item {[15]} A. Ivi\'c, The Laplace and Mellin transforms of powers of the Riemann
zeta-function, International Journal of Mathematics and Analysis
{\bf1}(2), 2006, 113-140.

\item {[16]} A. Ivi\'c, Remarks on the natural boundary of $\Z_k(s)$,
notes of October 2006 and October 2009.

\item {[17]} A. Ivi\'c, On some reasons for doubting the Riemann Hypothesis,
in ``The Riemann Hypothesis'', P. Borwein et al.,  CMS Books in Mathematics,
Springer, 2008, pp. 131-160.

\item{[18]} A. Ivi\'c, M. Jutila and Y. Motohashi,
The Mellin transform of powers of the zeta-function, Acta
Arithmetica {\bf95}(2000), 305-342.

\item {[19]}  M. Jutila, The Mellin transform of the square of Riemann's
zeta-function,  Periodica Math. Hungarica {\bf42}(2001), 179-190.

\item{[20]} M. Jutila, Atkinson's formula for Hardy's function,
J. Number Theory {\bf129}(2009), 2853-2878.

\item{[21]} A.A. Karatsuba and S.M. Voronin, The Riemann zeta-function,
Walter de Gruyter, Berlin etc., 1992.

\item{[22]} M.A. Korolev, On the integral of Hardy's function $ Z(t)$,
Izv. Math. {\bf}72, No. 3, 429-478 (2008); translation from Izv.
Ross. Akad. Nauk, Ser. Mat. {\bf72}, No. 3, 19-68 (2008).

\item{[23]} M. Lukkarinen, The Mellin transform of the square
of Riemann's zeta-function and Atkinson's formula, Ann.
Acad. Sci. Fenn. Math. Diss. {\bf140}, 2007.

\item {[24]} Y. Motohashi,  Spectral theory of the Riemann
zeta-function, Cambridge University Press, Cambridge, 1997.

\item {[25]} Y. Motohashi, The Riemann zeta-function and Hecke
congruence subgroups II, Journal of Research Institute
of Science and Technology, Tokyo, 2009, to appear.

\item{[26]} K. Ramachandra, On the mean-value and omega-theorems
for the Riemann zeta-function, Tata Inst. of Fundamental Research,
(distr. by Springer Verlag, Berlin etc.), Bombay, 1995.

\item {[27]} E.C. Titchmarsh, Introduction to the Theory of Fourier
Integrals,  Clarendon Press, Oxford, 1948.

\item{[28]} E.C. Titchmarsh, The theory of the Riemann
zeta-function (2nd edition),  University Press, Oxford, 1986.

\endRefs

\enddocument

\bye